\documentclass[11pt]{amsart}
\usepackage[all]{xy}
\usepackage{xspace,amssymb,amsfonts,epsfig,euscript,eufrak}
\author{Ben Webster}
\title[Khovanov-Rozansky homology via a canopolis formalism]{Khovanov-Rozansky homology\\ via a canopolis formalism}
\address{Department of Mathematics\\
         University of California, Berkeley\\
         Berkeley, CA 94720}
\subjclass[2000]{Primary 57M27; Secondary 13D02}
\email{bwebste@math.berkeley.edu}
\urladdr{http://math.berkeley.edu/~{}bwebste}
\thanks{This material is based upon
  work supported under a National Science Foundation Graduate Research
  Fellowship and partially supported by the RTG grant DMS-0354321.}
\begin{document}
  \newcommand{\nc}{\newcommand}
  \newcommand{\renc}{\renewcommand}
  \nc{\ti}{\tilde}
  \nc{\al}{\alpha}
  \nc{\mcC}[2]{\mathcal{C}^{#1}_{#2}}
  \nc{\ep}{\epsilon}
  \nc{\maf}[2]{\mathsf{MF}^{#1}_{#2}}
  \nc{\MF}{\mathbf{MF}}
  \nc{\KR}{\mathcal{R}}
  \nc{\KRh}{\mathcal{K}}
  \nc{\Kom}{\mathrm{Kom}}
  \nc{\mc}[1]{\mathcal{#1}}
  \nc{\Z}{\mathbb{Z}}
  \nc{\vp}{\varphi}
  \nc{\By}{\mathbf{y}}
  \nc{\Bx}{\mathbf{x}}
  \nc{\Bt}{\mathbf{t}}
  \nc{\No}{N_1}
  \nc{\Nt}{N_2}
  \nc{\Kmf}[2]{\ti{\mc Z}_{#1,#2}}
  \nc{\mf}{matrix factorization\xspace}
  \nc{\mfs}{matrix factorizations\xspace}
  \nc{\qpm}{near-isomorphism\xspace}
  \nc{\qpim}{near-isomorphism\xspace}
  \nc{\qpms}{near-isomorphisms\xspace}
  \nc{\dpl}{d_+}
  \nc{\dm}{d_-}
  \nc{\dpm}{d_\pm}
  \nc{\fr}[1]{\mathfrak{#1}}
  \nc{\Ht}{H^t}  
  \nc{\fH}{\EuScript{H}}
  \nc{\nai}{\EuScript{N}}
 \nc{\excise}[1]{}
 \nc{\Hhv}[2]{H_{hv}^{#1}\left(#2\right)}
 \nc{\Hv}[2]{H_{h}^{#1}\left(#2\right)}
 \nc{\hd}[1]{\EuScript{K}\left(#1\right)}
 \nc{\Span}{\mathrm{span}\,}
 \nc{\Hh}[2]{H_{v}^{#1}\left(#2\right)}
  \newtheorem{defi}{Definition}
  \newtheorem{thm}{Theorem}[section]
  \newtheorem{lem}[thm]{Lemma}
  \newtheorem{prop}[thm]{Proposition}
  \newtheorem{cor}[thm]{Corollary}
  
  \theoremstyle{remark}
  \newtheorem{remark}{Remark}
\begin{abstract}
  In this paper, we describe a canopolis (i.e. categorified planar
  algebra) formalism for Khovanov and Rozansky's link homology theory.
  We show how this allows us to organize simplifications in the matrix
  factorizations appearing in their theory.  In particular, it will
  put the equivalence of the original definition of Khovanov-Rozansky
  homology and the definition using Soergel bimodules in a more
  general context, allow us to give a new proof of the invariance of
  triply graded homology and give new analysis of the behavior of
  triply graded homology under the Reidemeister IIb move.
\end{abstract}

\maketitle

In the papers \cite{KR04,KR05}, Khovanov and Rozansky introduced a series
of homology theories for links.  These theories categorify the quantum
invariants for $\fr{sl}_n$, and the HOMFLYPT polynomial.  Unfortunately,
they remain very difficult to calculate, not least because
of the complicated matrix factorizations used in their original
combinatorial definition.  Later work of I. Frenkel, Khovanov, and Stroppel
\cite{FKS05,Kho05,Str06b,Str06a} has suggested a more systematic
definition of these invariants, and connection between these theories
and the structure of the BGG category $\mc O$ for the Lie algebra
$\fr{gl}_n$, but progress toward computational simplifications along
these lines has been slow.

In this paper, we will show that these invariants can be understood,
computed and in fact, defined in the context of canopolises.  We hope that
this approach will both lead to computational benefits and help the
reader to understand the definition of Khovanov-Rozansky homology
better.  A {\bf canopolis}\footnote{Bar-Natan uses the term ``canopoly''
  in the published version of \cite{BarN05}, but the consensus choice
  now seems to be the more etymologically correct ``canopolis.''} is a
categorification of the notion of a planar algebra defined by Bar-Natan
\cite{BarN05} (see Section~\ref{sec:canopolises}).  

Consider a disk in the plane with $m$ disks removed from its interior
(we call the places left by these removed disks ``holes'').  An {\bf
  oriented planar arc diagram} (or ``spaghetti-and-meatballs
diagram'') $\eta$ on this disk is a collection of oriented simple
curves with endpoints on the boundary of the disk (including the
boundary of the holes), along with choice of a distinguished point on
each boundary of a component (in diagrams, this point is distinguished
by putting a star next to it), and an ordering of the holes of the
diagram.

Let $\mc Q_i(\eta)$ be the set of planar arc diagrams $\omega$ such that
the outer boundary of $\omega$ matches the boundary of the $i$th hole of
$\eta$.  That is, there is the same number of endpoints, and if we order
the endpoints, starting at our distinguished point, the orientations of
the arcs match.

\begin{figure}[htbp]
  \centering
  \centerline{\epsfig{figure=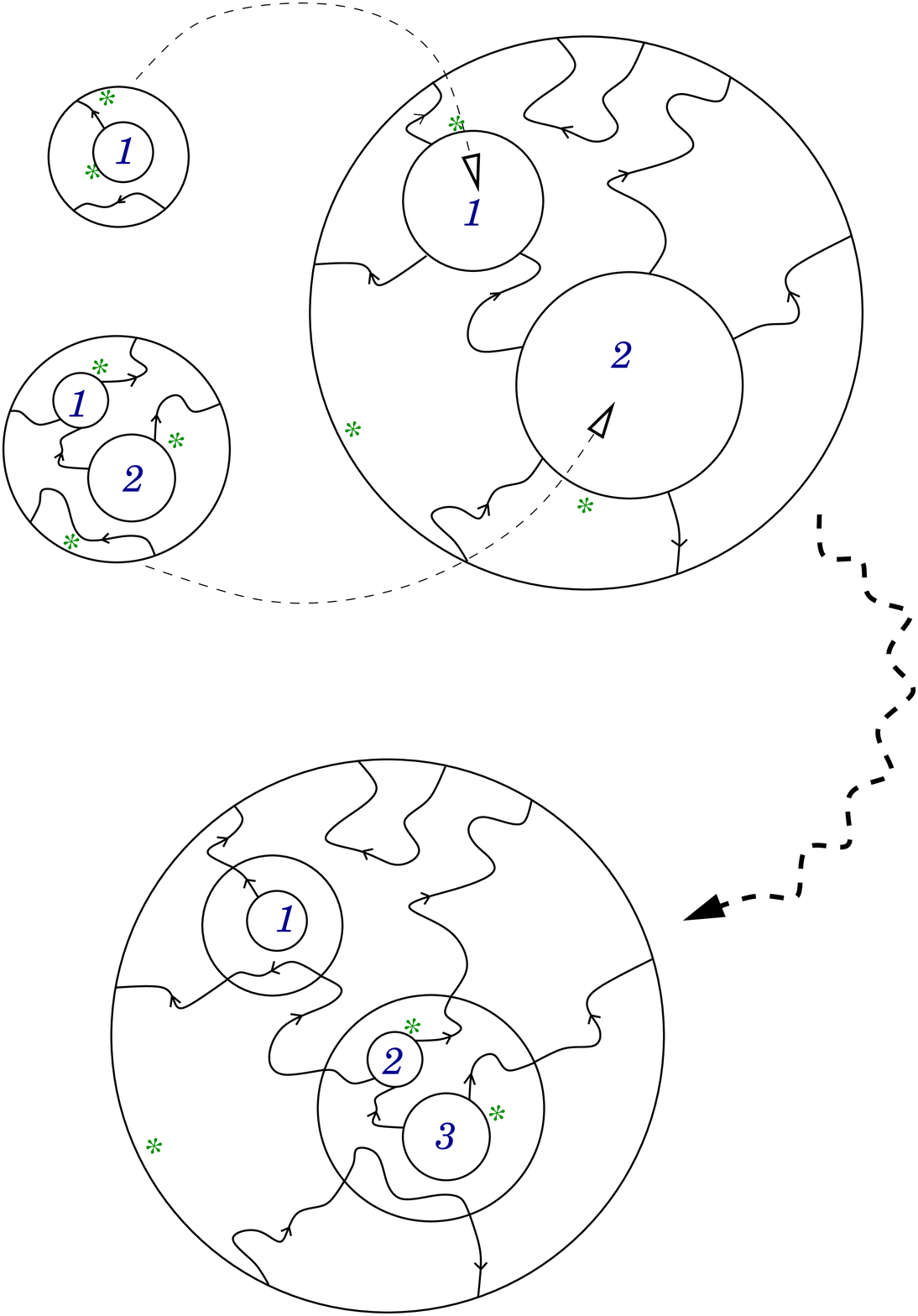, height=18cm}}
  \caption{The action of planar arc diagrams.}
  \label{fig:1}
\end{figure}

Any fixed planar arc diagram $\eta$ with a distinguished hole defines an
operation $\ti\eta: \mc Q_1(\eta)\times\cdots \times \mc Q_m(\eta)\to\mc
Q_0(\eta)$, by shrinking the given $m$-tuple of diagrams, and pasting it
into the holes of $\eta$, as is shown in Figure~\ref{fig:1}.  This
``multiplication'' is a particular instance of an algebraic structure
called an {\bf colored operad}. 

\begin{figure}[htbp]
  \centering
   \centerline{\epsfig{figure=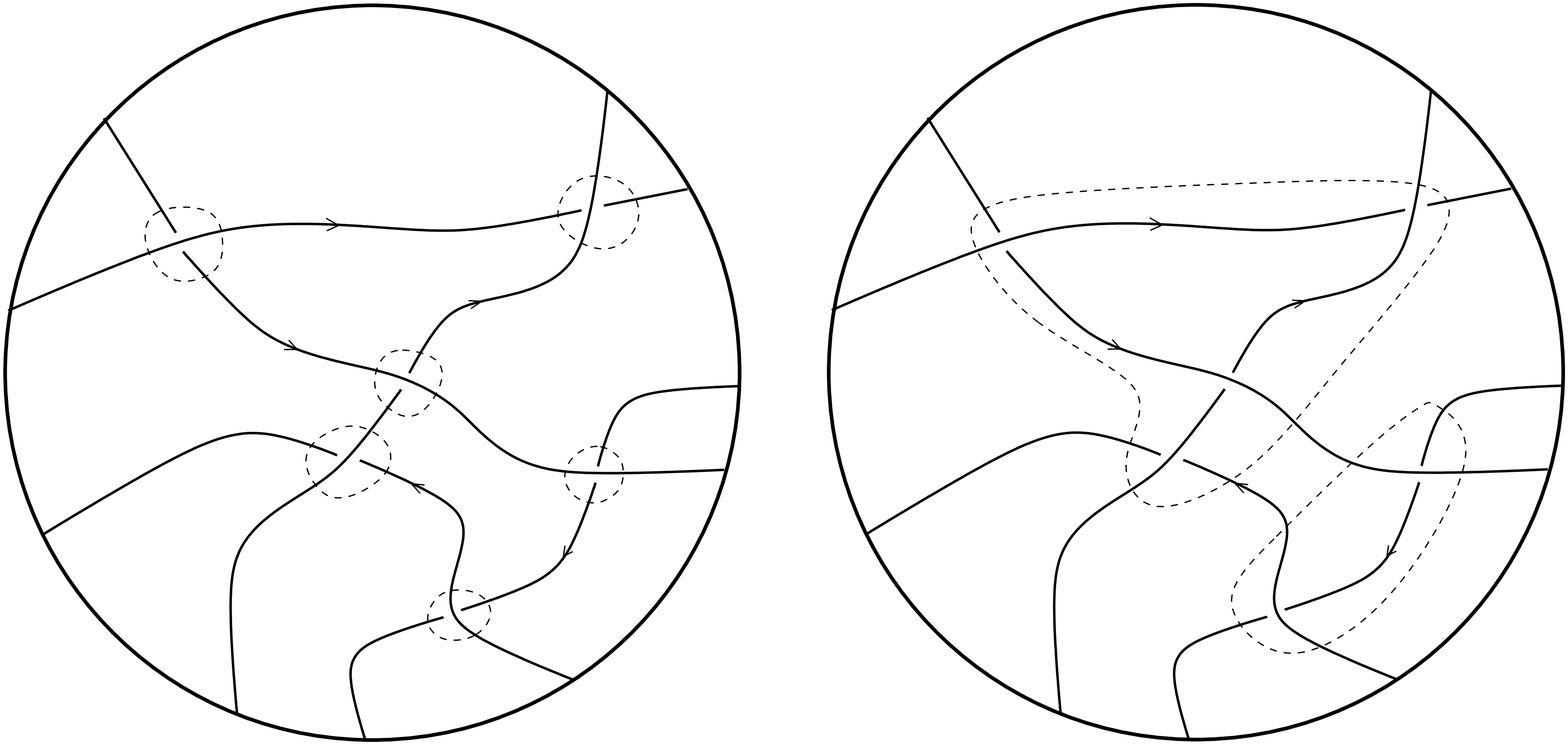, height=5.5cm}}
  \caption{Factoring a tangle (a) as a product of crossings and a planar arc
    diagram or (b) as a more general factorization.}
  \label{fig:2}
\end{figure}

The operad of planar arc diagrams acts on tangle diagrams in a disk as
well.  Phrased in the language of \cite{BarN05}, the set of tangle
diagrams $\EuScript{T}_{S,\ep}$ in a disk, with endpoints on the
boundary and a marked point on the boundary, partitioned according to
the orientation $\ep:S\to\{\pm 1\}$ of the endpoints, form {\bf a
  planar algebra}, that is, a set on which the operad of planar arc
diagrams acts.  In fact, one can build any tangle diagram from single
crossings in a disk and the action of a planar arc diagram.  More
generally, we will be interested in factoring a tangle as the action
of a planar arc diagram on simpler tangle diagrams. We depict these
operations in Figure~\ref{fig:2}.

One fruitful approach to the Jones polynomial and other quantum
invariants is to regard each as a homomorphism of planar algebras.
Thus, one can compute the Jones polynomial of a tangle once (recall
that this is an element of a certain vector space over $\mathbb{C}(q)$, rather
than just a polynomial), and then whenever one wishes to the Jones
polynomial of a knot, one cuts it into a planar arc diagram acting on
tangles whose Jones polynomial is known.

This approach is of more than theoretical value; if programed
skillfully, it can be extremely efficient.  Bar-Natan
\cite{BarN02,BarN05} presented a beautiful extension of this approach
to Khovanov's original link homology (the $\fr{sl}_2$-version of
Khovanov-Rozansky), which at once gives a simple description of the
knot homology and is extremely computationally efficient, allowing the
computation by computer of Khovanov homology for knots of dozens of
crossings.

Unfortunately, we do not know how to give a similar, matrix
factoriza\-tion-free description of Khovanov-Rozansky. Instead, we will
show that Khovanov-Rozansky can be defined using certain homotopy
categories of matrix factorizations which admit an action of planar
arc diagrams, that is, a canopolis structure, which we can see as an
analogue of Bar-Natan's geometric canopolis.

While this is an essentially formal construction, it allows us to
simplify the KR-complex of a small tangle before or after we apply the
action of a planar arc diagram (where ``simplify'' has a very precise
definition, given in Section~\ref{sec:simplifications-krt}), allowing us
to organize computations of KR homology according to Bar-Natan's
``divide and conquer'' philosophy. In particular, it will give us a new
understanding of the equivalence between KR homology and the homology
defined by Soergel bimodules, shown by Khovanov in the HOMFLYPT case
\cite{Kho05}.

We will also apply this approach to the triply graded homology theory
discussed in \cite{Kho05, KR05} to give a new proof of invariance and
show that the changes in triply graded homology when the diagram undergoes
a second Reidemeister move is controlled by a certain spectral
sequence.

\section{Matrix Factorizations}

\subsection{Preliminaries on \mfs}

We will attempt to follow the notations and terminology of Rasmussen
\cite{Ras06}.  Let $M$ be a $\Z$-graded module over a ring $\mc S$.

\begin{defi}
A \emph{(}$\Z$-graded\emph{)} {\bf matrix factorization} on $M$ with
potential $\vp\in \mc S$ is a map $d=\dpl+\dm:M\to M$ with $\dpm$ of
graded degree $\pm 1$ such that $d^2=\vp$.
\end{defi}

Though this is not usual definition of a matrix factorization (where
typically we only assume a $\Z/2$ grading), this richer structure is
more useful from the perspective of knot theory. 

Matrix factorizations over $\mc S$ with a fixed potential naturally form an
abelian category, with morphisms given by maps commuting with $d$.  We
only assume that these maps are homogeneous with respect to the
$\Z/2$-grading.  

Even better, we
can think of all matrix factorizations over all unital rings as a
2-category $\MF$ such that:
\begin{itemize}
\item The objects are given by a pairs of a unital ring, and an element
  of that ring $(\mc S,\vp)$.
\item The 1-morphisms between $(\mc S_1,\vp_1)$ and $(\mc S_2,\vp_2)$
  given by matrix factorizations over $\mc S_1\otimes\mc  S_2$ with
  potential $\vp_1\otimes 1-1\otimes \vp_2$, and composition given
  by tensor product of matrix factorizations of bimodules.
\item The 2-morphisms between two matrix factorizations given by morphisms
  in the usual sense.
\end{itemize}
While matrix factorizations may seem strange, in fact, they arise
very naturally in homological algebra (see, for example
\cite{DE80}). Consider a module $N$, and a ring element $\vp\in\mc
S$ which annihilates $N$.  Fix a finite-length free resolution
$\mathbf{M}^\bullet$ of $N$ (by convention, we give this resolution
the cohomological grading, i.e. the differential is of degree 1),
and let $M$ be the direct sum of all its components. Let $\dpl:M\to
M$ be the differential, and $\vp_M$ be the action of multiplication
by $\vp$ on $M$. Since the induced map $\vp_N$ is 0, by standard
homological algebra, $\vp_M$ is homotopic to 0, i.e. there exists a
map $\dm:\mathbf{M}^\bullet\to\mathbf{M}^\bullet$ of degree 1 such
that $\dpl\dm+\dm\dpl=\vp$.

Now, assume that $\dm$ also defines the structure of a chain complex on
$\mathbf{M}^{\bullet}$, that is, $\dm^2=0$. Let $d=\dpl+\dm$.
By the homotopy formula above, we have 
\begin{equation*}
  d^2= \dpl^2+\dpl\dm+\dm\dpl+\dm^2=\vp,
\end{equation*}
that is, $d$ defines a matrix factorization with potential $\vp$ on $M$
with the grading given by homological degree.

Recall for an ordered $n$-tuple $(x_1,\ldots,x_n)\in R^n$ in a
commutative ring $R$ is called a {\bf regular sequence} if the action
of $x_i$ is a non-zerodivisor (multipication by it is injective) on $R/(x_1,\ldots,x_{i-1})$ for all $i$.

In the case when $N$ is the quotient $\mc S/({\Bx})$ of $\mc S$ by the
ideal generated by a regular sequence $\Bx=\{x_1,\ldots,x_m\}$, the
matrix factorizations constructed from $N$ have a particularly nice
interpretation, as described in slightly different language by Eisenbud
in \cite[\S 17.4]{Ei97}. Let $Z_{x_i}$ denote the two-term complex $\mc
S\overset{x_i}{\longrightarrow}\mc S$. If $\mc S$ is graded and each
$x_i$ is homogeneous, then we can shift gradings so that this
differential is of degree 2 (while this may seem like a peculiar choice,
it the one which makes this grading match with the standard variable $q$
in the quantum invariants which KR homology should categorify).  This is
clearly a free resolution of $\mc S/(x_i)$.  Thus, the easiest possible
guess for a free resolution of $N$ is the Koszul complex
\begin{equation*}
   \mc Z_{\Bx}= \bigotimes_{i=1}^{m}Z_{x_i}.
\end{equation*}
This complex is indeed a resolution of $N$, since $\Bx$ is regular
\cite[\S 17.2]{Ei97}
(otherwise, it might have higher cohomology).

Since $\vp_N=0$, we must have $\vp=\sum_{i=1}^my_ix_i$ for some
sequence (not necessarily unique or regular) $\By=\{y_1,\ldots,y_n\}$. Instead
of taking the complex $Z_{x_i}$, consider the matrix factorization
\begin{equation*}
  \ti Z_{x_i,y_i}=\xymatrix{\mc S \ar@/^/[r]^{x_i}&\mc S\ar@/^/[l]^{y_i}}.
\end{equation*}
By analogy with the Koszul resolution, we define the {\bf Koszul matrix
factorization} of the pair $(\Bx,\By)$ to be the tensor product
\begin{equation*}
\Kmf\Bx\By=\bigotimes_{i=1}^{m}\ti Z_{x_i,y_i}.
\end{equation*}
This is a matrix factorization with potential $\sum_ix_iy_i=\vp$.

\subsection{Near-isomorphisms}

We would like to generalize the following standard fact of homological
algebra:
\begin{prop}\label{T:tens-qi}
  Let $f:C\to C'$ be a chain map between complexes of $\mc  S$ modules
  which induces an isomorphism on homology, that is, a
  quasi-isomorphism. Then for any complex of $D$ of projective $\mc
  S$-modules, $f\otimes 1: C\otimes D\to C'\otimes D$ is also a
  quasi-isomorphism.
\end{prop}
\begin{proof}
  We use a spectral sequence on $C\otimes D$, associated to the
  filtration
  \begin{equation*}
    F_n=\bigoplus_{i\leq n} C_j\otimes D_i
  \end{equation*}
  In this case $F_{n}/F_{n-1}\cong C\otimes D_n$.  Since $f\otimes 1$ is
  filtered, it induces a map of spectral sequence, and since $D_n$ is
  projective, the induced map on the $E^0$ term is an isomorphism.  Thus
  it is an isomorphism on the $E^\infty$ term as well, which is the
  associated graded module of $H^*(C\otimes D)$.  Thus, $f$ is a
  quasi-isomorphism.
\end{proof}

Note that this result depends heavily on the fact that $D$ is
projective (or more generally, flat).

How might such a fact be generalized to matrix factorizations?  First,
let us define a class of maps analogous to quasi-isomorphisms of
complexes (unfortunately ``quasi-isomorphism of matrix factorizations,''
as defined by Rasmussen \cite{Ras06} means something slightly different,
and more analogous to a homotopy equivalence of complexes).

\begin{defi}
  We call a chain map $f:M_+\to M'_+$ a \qpm if for each \mf $D$ on a
  projective $\mc S$-module of potential $-\vp$, the map
  $\Ht(M\otimes D)\to\Ht(M'\otimes D)$ is an isomorphism.
\end{defi}

Unlike in the case where $\vp=0$, the hard part now will be identifying such
morphisms.  While we know of no explicit characterization, one example
will be sufficient for our purposes

Let $M$ be a matrix factorization of potential $\vp$ with $M^i=0$ for
all $i>0$, and let $\fH(M)$ be $H^0(M_+)$.  Obviously, there is a
chain map $\pi:M_+\to \fH(M)$, where $\fH(M)$ is considered as a complex
concentrated is degree 0. Note that $\fH(M)$ must be annihilated by
$\vp$, since the image of $\dpl$ in $M^0$ contains $\vp\cdot M^0$.

\begin{thm}\label{res-qpm}
  If the natural map $\pi:M_+\to \fH(M)$ is a quasi-isomorphism of
  complexes, then it is also a \qpm of \mfs.
\end{thm}
\begin{proof}
  Let $K=\ker \dpl\subset M$.  Fix a matrix factorization of projective
  modules $D$ of potential $-\vp$, and consider the following filtration
  on $M\otimes D$:
  \begin{equation*}
    G^n=\left(\bigoplus_{i<n, j\in\Z} M^i\otimes D^j\right)\oplus
    \left(\bigoplus_{j\in\Z}K^n\otimes D^j\right)
  \end{equation*}
  This filtration is a mix of the standard choice on a tensor product of
  complexes, and that used by Rasmussen for any matrix factorization of
  potential zero \cite[Lemma 5.11]{Ras06}.

  Since $d=\dpl+\dm$ preserves this filtration, there is an associated
  spectral sequence, converging to $\Ht(M\otimes D)$.

  The $E^0$ term of this spectral sequence is 
  \begin{equation*}
    G_{n+1}/G_{n}\cong \left(K^{n+1}\oplus M^n/K^n\right)\otimes D
  \end{equation*}
  and as Rasmussen computed, the differential $d_0$ is given by the
  total differential on the above tensor product, when we make
  $K^{n+1}\oplus M^n/K^n$ into a matrix factorization by 
  \begin{equation*}
    \xymatrix{K^{n+1}\ar@/^/[r]^{\dm^M}&M^n/K^n\ar@/^/[l]^{\dpl^M}}.
  \end{equation*}
  Since we assumed that $H^i(M_+)=0$ for all $i\neq 0$, the map $\dpl^M$
  induces an isomorphism $M^i/K^i\cong K^{i+1}$ for all $n\neq 0$, the
  complex $\left(G_{n+1}/G_{n}\right)_+$ is exact and the
  homology of $ G_{n+1}/G_{n}$ is trivial if $n\neq 0$.

  Thus, all higher differentials have trivial source or trivial target,
  so our sequence collapses at $E^1$.  Since $G^1/G^0\cong \fH(M)\otimes D$,
  we find that 
  \begin{equation*}
    H^t(M\otimes D)\cong H^t(\fH(M)\otimes D),
  \end{equation*}
with the
  isomorphism induced by the natural projection $\pi$.
\end{proof}

\subsection{Khovanov-Rozansky \mfs}\label{sec:khov-rozansky-mfs}

Matrix factorizations appear in knot theory through the work of Khovanov
and Rozansky: they associate to any oriented tangle diagram $T$ and any
polynomial $p$, which vanishes at 0, a complex of matrix factorizations
we denote by $\KR_p(T)$, defined as follows:

Consider the graph $\mc G(T)$ of $T$ which has vertices corresponding to
crossings or end points of components of $T$ and edges corresponding to
segments of diagram between crossings. The orientation of $T$ induces an
orientation on $\mc G(T)$. Let $\mc F(T)$ denote the set of flags of the
graph $\mc G(T)$, that is, pairs $(x,e)$ of adjacent edges and vertices.
Let $\mc{\ti S_T}=\ti{\mc S}$ be polynomials over $k$ in the variables
$t_{x,e}$, where $(x,e)$ ranges over $\mc F(T)$.

\begin{figure}[htbp]
  \centering
    \centerline{\epsfig{figure=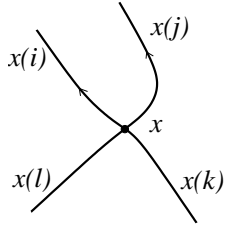, height=3cm}}
  \caption{The labeling of edges around a crossing}
  \label{fig:3}
\end{figure}

For each crossing $x$, number the adjacent edges as shown in
Figure~\ref{fig:3}. Let $t_h=t_{x,e_h}$, where $h=i,j,k,\ell$. 

Let $L_x'$ be a Koszul \mf on the Koszul resolution of
\begin{equation*}
  \Bt_x'=(t_i+t_j-t_k-t_\ell,t_it_j-t_kt_\ell)
\end{equation*}
with potential $\vp_x=p(t_i)+p(t_j)-p(t_k)-p(t_\ell)$.  Recall that we
have fixed the grading on Koszul resolutions such that $\dpl$ has graded
degree 2 and the $0$th term has the same grading as the ring itself with
no shift.  Such a factorization exists, since $p(t_i)+p(t_j)$ is a
polynomial in $t_i+t_j$ and $t_it_j$, by the fundamental theorem of
symmetric function theory.  Its exact form will not be important for us
at the moment.  Since $\Bt_x'$ is regular, $(L_x')_+$ is a free
resolution of $\fH(M)\cong \mc {\ti S}/(\Bt_x')$, where $(\Bt_x')$ denotes the
ideal in $\mc {\ti S}$ generated by the elements of $\Bt_x'$.

Let $L''_x$ be a Koszul matrix factorization with potential $\vp_x$ on
the Koszul complex of $\Bt''_x=(t_i-t_\ell,t_j-t_k)$.  As with $\Bt'_x$,
this is a regular sequence, so $(L_x'')_+$ is a free resolution of
$\fH(M)\cong \mc {\ti S}/(\Bt''_x)$.

Khovanov and Rozansky define a two term complex $L_x$ which depends on whether
the crossing was positive or negative;  if $x$ is positive, it is a
complex of the form $\rho^+_x:L'_x\to L_x''$ and if $x$ is negative, it is of
the form $\rho^-_x:L_x''\to L_x'\{-2\}$, where $\{a\}$ denotes grading
shift by $a$, in each case with $L_x''$ in
homological degree 0.  The exact form of these maps will not be
important to us at the moment.  We only note that after applying
$\fH(-)$, 
\begin{itemize}
\item the induced map $\fH(\rho^+_x):\mc {\ti S}/(\Bt')\to \mc {\ti S}/(\Bt'')$ is the
obvious projection.
\item the induced map $\fH(\rho^-_x):\mc {\ti S}/(\Bt'')\to \mc {\ti S}/(\Bt')\{-2\}$ is
  that induced by multiplication by $t_i-t_k$. 
\end{itemize}
It's worth noting that the complex $\fH(L_x)$ is independent of $p$.

For each edge $e$, directed from $x_a$ to $x_b$, we let $L_e$ be the
Koszul matrix factorization of the pair $\Bx=(t_a-t_{b})$ and
$\By=\left(\frac{p(t_{a})-p(t_{b})}{t_{a}-t_{b}}\right)$, where
$t_a=t_{x_a,e}$ and $t_b=t_{x_b,e}$, again with polynomial grading such
that $\dpl$ is of degree 2.

We define the Khovanov-Rozansky complex of the diagram $T$ to be the
complex of matrix factorizations given by the tensor product
\begin{equation*}
  \KR_p(T)=\left(\bigotimes_{e}L_e\right)\otimes \left(\bigotimes_{x}L_x\right).
\end{equation*}

Since each \mf in this complex is a direct sum of Koszul \mfs, we can
apply the functor $\fH$ component-wise, to obtain a complex of modules
over $\mc {\ti S}$, which we will refer to
as the ``naive Khovanov-Rozansky complex'' $\nai(T)=\fH(\KR_p(T))$.  As
with a single crossing, $\nai(T)$ is independent of $p$, since
$\KR_p(T)_+$ is independent of $p$.

Let $\mc L(T)\subset \mc F(T)$ be the set of flags containing vertices
of degree 1.  There is only one such flag for each vertex of degree 1,
and these correspond to the open ends of the tangle $T$.  Let $\mc S_T$
denote the subring of $\ti{\mc S}$ generated by the variables
corresponding to elements of $\mc L(T)$. Define $\ep:\mc F(T)\to\{\pm1\}$
by
\begin{align*}
  \ep(x,e)=
  \begin{cases}
    1 & e\text{ is directed into }x\\
    -1 & e\text{ is directed out of }x
  \end{cases}
\end{align*}

Then the potential of $\KR_p(T)$ is given by the sum
\begin{equation*}
    \vp_T=\sum_{(\ell,e)\in\mc L(T)}\ep(\ell,e)p(t_{\ell,e})
\end{equation*}
In particular, we can consider $\KR_p(T)$ as a complex of matrix
factorizations over $\mc S_T$, a change which seems small, but will be
key to simplifications we do later.

If $T$ is the diagram of a link (i.e. a closed tangle), then $\KR_p(T)$
has potential 0, and we can take the total homology of each matrix
factorization. In this case, we obtain a complex of graded vector
spaces, which we call $\KRh_p(T)$.  The homology of this complex (as a
bigraded vector space) is what is typically called {\bf unreduced
  Khovanov-Rozansky homology}.  We can obtain {\bf reduced
  Khovanov-Rozansky} (for a knot) by quotienting out by the action of
one of the generators of $\mc S_T$ on $\KRh_p(T)$, and then taking
homology.  If $p$ is homogeneous (i.e. $p(x)=x^{N+1}$ for some $N$) then
on both these homologies two gradings will survive (otherwise, we will
only have the homological grading).  We will not go into the details, as
they are not of great importance to the rest of the paper, and are
covered in great detail in \cite[\S 2]{Ras06}.

\section{Canopolises}
\label{sec:canopolises}

\subsection{Canopolises of matrix factorizations}
\label{sec:canopo-whats}

The theory of planar algebras originated with Vaughan Jones's theory of
subfactors \cite{Jon99}, and they have shown themselves to be a very
useful formalism for dealing with knot invariants.  In his reformulation
of Khovanov homology, Bar-Natan \cite{BarN05} uses a categorification of
a planar algebra, called a canopolis:

\begin{defi}
A (oriented) {\bf canopolis} is an assignment of
\begin{enumerate}
\item A category $\mcC{X}\ep$ for each set totally ordered finite set
  $X$ equipped with sign map $\ep:X\to \{+,-\}$.  We think of this as being
  associated to a disk with signed marked points on the boundary, and
  with a distinguished marked point (so points are totally ordered,
  not just cyclically ordered).
\item A functor
  $\ti\eta:\mcC{X_1}{\ep(1)}\times\cdots\times\mcC{X_m}{\ep(m)}
  \to\mcC{X_0}{\ep(0)}$, for each oriented planar arc diagram $\eta$,
  where $X_j$ denotes the set of endpoints of arcs on the $j$th boundary
  component, with the sign determined by the orientation of the arc at
  that point.  The action of planar arc diagrams should commute with
  composition of functors.
\end{enumerate}
\end{defi}

Bar-Natan's original examples were for the most part very geometric,
being modifications of various cobordism categories.  Rather than attempt
to do justice to his presentation, we refer the interested reader to his
paper \cite[pg. 31]{BarN05}. 

The geometric canopolis of interest to us will be as follows:
\begin{itemize}
\item We let $\mcC{X}\ep$ be the category of oriented tangles in a
  thickening of the disk $D^2$ to a 3-ball $\Sigma D^2=B^3$, with the
  endpoints of the tangles on the distinguished points $X$, and the
  orientation of the tangle matching the sign sequence $\ep$.
  Morphisms between $T$ and $T'$ are oriented cobordisms embedded in
  $I\times B_3$, with boundary given by $T\times\{0\}\cup
  T'\times\{1\}\cup X\times I$.
\end{itemize}
We denote this canopolis by $\EuScript{C}$

As promised, for any polynomial $p$ (vanishing at 0, as before), we will
define an associated canopolis $\EuScript{M}_p$ of \mfs, which is a
natural home for KR homology.

First, associated to the sign sequence $\ep:X\to \{+,-\}$ is the
category $\maf{X,p}{\ep}$ of matrix factorizations over a polynomial
ring $k[X]$ with generators indexed by $C$ of potential $\sum_{x\in X}
\ep_xp(t_x)$.  

As we discussed before, the most natural functors from
$\maf{X_1,p}{\ep(1)}\times \cdots \times\maf{X_m,p}{\ep(m)}$ to
$\maf{X_0,p}{\ep(0)}$ are those induced by tensor product with a matrix
factorization over $k[X_0,\ldots, X_m]$ with potential
\begin{equation*}
  \sum_{x\in X_0}\ep_x(0)p(t_x)-\sum_{j=1}^m\sum_{x\in X_j}\ep_x(j)p(t_x)
\end{equation*}

In fact, there is a clear choice in this category: Let $\mc A(\eta)$ be
the set of arcs of $\eta$.  Each arc $\al\in \mc A(\eta)$ has a head
$\al_+$ and a tail $\al_-$.  

Now define sequences $\Bx,\By$ by 
\begin{align*}
  \Bx&= \big(t_{\al_+}-t_{\al_-}\big)_{\al\in \mc
    A(\eta)}\\
  \By&=\left(\frac{p(t_{\al_+})-p(t_{\al_-})}
    {t_{\al_+}-t_{\al_-}}\right)_{\al\in \mc A(\eta)}
\end{align*}
and let $\ti{\mc Z}_\eta$ be the Koszul matrix factorization
of this pair.  This is a matrix factorization over $k[X_0,\ldots, X_m]$,
and its potential is
\begin{equation*}
  \sum_{\al\in\mc A(\eta)}p(t_{\al_+})-p(t_{\al_-})= \sum_{x\in X_0}\ep_x(0)p(t_x)-\sum_{j=1}^m\sum_{x\in X_j}\ep_x(j)p(t_x)
\end{equation*}
since each marked point on any boundary of $\eta$ is the endpoint of
exactly one arc.

The canopolis functor $\ti\eta: \maf{X_1,p}{\ep(1)}\times \cdots
\times\maf{X_m,p}{\ep(m)} \to\maf{X_0,p}{\ep(0)}$ will simply be tensor
product with $\ti{\mc Z}_\eta$ over $k[X_1,\ldots, X_m]$ .

Note that this also induces a canopolis structure on the categories of complexes of matrix factorizations $ \mathrm{Kom}(\maf{X,p}{\ep})$ and the homotopy category of complexes $\hd{\maf{X,p}{\ep}}$, since tensoring with $\ti{\mc Z}_\eta$ is exact.

\begin{thm}
  These functors define a canopolis structure $\EuScript{M}_p$ on
  $\hd{\maf{X,p}{\ep}}$.  Furthermore,
  $\KR_p:\EuScript{C}\to\EuScript{M}_p$ is a functor of canopolises,
  i.e. the diagram
\begin{equation*}
  \xymatrix{\displaystyle{\prod_{j=1}^m \EuScript{T}_{X_j,\ep(j)}}
    \ar[d]_{\KR_p}\ar[r]^{\ti\eta}
    & \EuScript{T}_{X_0,\ep(0)}\ar[d]^{\KR_p}\\
    \displaystyle{\prod_{j=1}^m
    \hd{\maf{X_j,p}{\ep(j)}}} 
  \ar[r]_{\ti\eta}
  &\hd{\maf{X_0,p}{\ep(0)}}}
\end{equation*}
is commutative.  In particular, if $\ti\eta(T)$
is closed, then $\KRh_p(\ti\eta(T))$ and $H^*(\ti\eta_p(\KR_p(T)))$ are
isomorphic as complexes.
\end{thm}
Note that while $\KR_p$ is a $\Z$-graded matrix factorization, the maps
associated to cobordisms typically only preserve the $\Z/2$-grading.
\begin{proof}
  Luckily, all the necessary computations were done by Khovanov and
  Rozansky in \cite{KR04}.  Checking that the composition of planar arc
  diagrams matches with composition of functors is simply applying 
  \cite[Prop. 15]{KR04} at each pair of boundary points which are glued
  together. 
  
  The commutation with the functor $\KR_p$ is simply rephrasing the
  original definition, after placing a mark on each connected pair of
  boundary points.
\end{proof}

In particular, though $\KR_p(T)$ was first defined over a ring 
$\mc{\ti S}_T$ with variables correspond to all elements of $\mc F(T)$,
this canopolis formalism shows that we need only remember the action of
variables corresponding to endpoints, not to internal edges of $\mc G(T)$.
Often after restricting to this smaller subring, we can identify trivial
summands of the complex $\KR_p(T)$.

\subsection{Simplifications in $\KR_p(T)$}
\label{sec:simplifications-krt}

For our purposes, a simplification of a \mf will be a quotient such that
the projection map is a \qpm.  

Consider a Koszul matrix factorization $M=\mc{\ti Z}_{\Bx,\By}$.  Then
we expect $M$ to have a great number of simplifications.  For any $n<m$,
we can rewrite $M$ as a tensor product $M\cong M'\otimes_{\mc S} M''$, where
$M'={\mc{\ti Z}_{\Bx',\By'}}$ and $M''=\mc{\ti Z}_{\Bx'',\By''}$, and
\begin{align*}
  \Bx'&=(x_{1},\ldots,x_{n}), & \By'&=(y_{1},\ldots,y_{n}),\\
\Bx''&=(x_{n+1},\ldots, x_{m}), & \By''&=(y_{n+1},\ldots,y_{m}).
\end{align*}

Now, assume $\Bx'$ is a regular sequence.  In this case $M'_+$ is a
free resolution of $\fH(M')=\mc S/(\Bx')$.  By Theorem~\ref{res-qpm},
the natural map $\pi:M_+'\to \fH(M')$ is an \qpm, and thus,
$\pi\otimes 1:M\to \fH(M')\otimes_{\mc S}M''$ is as well.

This is a very useful principle in Khovanov-Rozansky homology.  For
instance, it gives a new proof of the equivalence of Rasmussen's
definition of KR homology with the original definition:  simply
apply the above construction the subsequence
$(t_{e,\al(e)}-t_{e,\omega(e)})_{E\in\mc G(T)}$, which appears in the
sequence for each term of the Khovanov-Rozansky complex.

We will concentrate on the dual approach of simplifying the \mfs
corresponding to crossings.  This explains why we want a different notion
of equivalence for \mfs from Rasmussen's: his approach was adapted to
keeping projective matrix factorizations on crossings, and tensoring
them with non-projective modules on edges, whereas ours is adapted to
having non-projective complexes on crossings, and projective \mfs on
edges.

First of all, we note that simplifications are preserved by the
canopolis action.

\begin{prop}\label{canop-qpm}
  If $f$ is a \qpm, then $\ti\eta(f)$ is also a \qpm.  Thus if $f$ is a
  simplification, so is $\ti\eta(f)$.
\end{prop}
\begin{proof}
  For the first statement, note that $\ti{\mc Z}_{\eta}$ is a \mf on a
  projective module, and thus tensor product with it preserves \qpms.
  For the second, we need only recall that tensor product is right
  exact.
\end{proof}

Thus, if we would like to calculate $\KRh_p(T)$ for some link diagram $T$,
but do not know how to simplify $\KR_p(T)$, then we might hope to factor
$T$ as $T=\ti\eta(T'_1,\ldots,T'_m)$, where the $T_j'$ are simpler
tangles for which we \emph{can} simplify $\KR_p(T)$, and then apply the
action of our canopolis.

For example, we have a simplification of $\KR_p(T_{\pm})$, where $T_\pm$ is a
disk with a single positive or negative crossing.  In this case,
$\KR_p(T_\pm)$ is simply the two term complex $L_x$ corresponding to the single
vertex in the graph of its projection.

\begin{prop}\label{crossing-qpm}
  The map $L_x\to\fH(L_x)$ is a degree-wise \qpm.
\end{prop}
\begin{proof}
  We noted in Section \ref{sec:khov-rozansky-mfs} that $L'_x$ and
  $L_x''$ are both \mfs on free resolutions, and the map
  $L_x\to\fH(L_x)$ is just the map to $\fH$ applied degree-wise (in
  the homological grading).  Thus, by Theorem \ref{res-qpm}, it is a
  degree-wise \qpm.
\end{proof}

If we factor a tangle into disks with a single crossing and a planar arc
diagram $\eta$, as shown in Figure~\ref{fig:2}(a), then $\KR_p(T)\cong
\ti\eta(\{L_x\})$ where $x$ ranges over crossings, ordered according to
the order chosen on the holes of $\eta$.

Let $\widehat{\KR}(T)$ be the complex $\ti\eta(\{\fH(L_x)\})\cong
\big( \bigotimes_{x}\fH(L_x)\big)\otimes \big(\bigotimes_e L_e\big)$.
Combining Propositions~\ref{canop-qpm}~and~\ref{crossing-qpm}, we find
that 
\begin{prop}\label{crossing-simp}
  The natural map $\KR_p(T)\to\widehat{\KR}_p(T)$ is a degree-wise \qpm
\end{prop}
Thus, we have a new complex of matrix factorizations which is the tensor
product of a complex of modules and a single regular Koszul matrix
factorization, and whose homology is $\KRh_p(T)$.

Of course, we may hope that we can simplify more general tangles than
single crossings using this philosophy.  For each tangle, we have a map
$\pi_T:\KR_p(T)\to \nai(T)$ from the honest KR complex to the naive one.
In Khovanov-Rozansky homology, as in life, things would be easier if we
could just be naive, but if we aren't savvy often enough, we can run
into trouble.  While Rasmussen's results show we can be ``savvy'' about
crossings and ``naive'' about edges, and Proposition~\ref{crossing-simp}
shows we can be ``naive'' about crossings and ``savvy'' about edges, we
will lose too much information if try take naive Khovanov-Rozansky of an
entire knot. After all, the naive Khovanov-Rozansky homology does not
depend on the polynomial $p$ and we know that honest Khovanov-Rozansky
homology does, since different choices of $p$ categorify
$\fr{sl}_n$ invariants for all $n$.

Thus, we would like to find a class of tangles about which we can be
naive, while still recovering honest KR homology.

\subsection{Acyclic tangles}
\label{sec:acyclic-tangles}

\begin{defi}
  We call an oriented tangle diagram {\bf acyclic} if the graph $\mc
  G(T)$ has no oriented cycles.
\end{defi}
Obviously a single crossing is acyclic.  Also, the tangles inside of the
dashed circles in  Figure~\ref{fig:2}(b) are acyclic, whereas the entire
tangle is not. Note that a tangle with a closed component is never
acyclic, whereas a braid always is.
 
\begin{thm}\label{acyc-qpm}
  If $T$ is acyclic, then $\pi_T:\KR_p(T)\to\nai(T)$ is an \qpm.
\end{thm}
\begin{proof}
  Since it is irrelevant to question at hand, we forget about the
  chain complex structure on $\KR_p(T)$ and consider it simply as a matrix
  factorization.  By Theorem~\ref{res-qpm}, we need only show that
  $H_i(\KR_p(T)_+)=0$, for any $i\neq 0$.
  
  We induct on the number of crossings.  Since $T$ is acyclic, the
  direction of edges induces a partial ordering on vertices.  Take any
  maximal element $x$. Using the conventions shown in
  Figure~\ref{fig:3}, $i(x)$ and $j(x)$ must be leaves or the adjacent
  vertex would be higher than $x$ in our partial order.  We will assume
  for simplicity that $k(x)$ and $\ell(x)$ are not leaves.  The case
  where they are follows from the same arguments we present below.
  
  Let $T'$ be $T$ with the crossing $x$ removed.  Thus, $\mc S_T\cong
  \mc S_{T'}\otimes k[x_i,x_j,x_k,x_\ell]$ where, as before, we let
  $t_h=t_{x,h(x)}$ where $h=i,j,k,\ell$, and
  \begin{equation*}
    \KR_p(T)_+\cong\KR_p(T')_+\otimes_k
  (L_x)_+\otimes_k (L_{\ell(x)})_+\otimes_k (L_{k(x})_+
  \end{equation*}
  Since $\KR_p(T')_+\otimes_k (L_x)_+$ is projective as a $\mc S_T$ module,
  we can project $(L_{\ell(x)})_+$ and $(L_{k(x})_+$ to their cohomology
  without changing the cohomology of $\KR_p(T)$.  This simplification is
  isomorphic to $\KR_p(T')_+\otimes_{k[t_\ell,t_k]} (L_x)_+$ where we
  let $k[t_\ell,t_k]$ act on $\KR_p(T)_+$ by the variable
  $t_{\ell'},t_{k'}$ corresponding to the other end of the edges
  $\ell(x),k(x)$.
  
  Since both $\KR_p(T')_+$ and $(L_x)_+$ are projective resolutions, the
  cohomology of this complex is
  $\mathrm{Tor}^i_{k[x_\ell,t_k]}(\nai(T),\fH(L_x))$.
  
  Both $\fH(L_x')$ and $\fH(L_x'')$ are free as $k[x_\ell,t_k]$-modules,
  as was proved by Soergel \cite{Soe92}.  Thus, all higher
  $\mathrm{Tor}$'s vanish, and we are done.
\end{proof}

\begin{cor}
  Let $T$ be a tangle diagram which can be factored as the action of $\eta$ on a
  set of acyclic tangles $\{T_i\}$, then the natural map $\KR_p(T)\cong
  \ti\eta(\{\KR_p (T_i)\}) \to \ti\eta(\{\nai(T_i)\})$ is a \qpm.  If
  $T$ is a link diagram, then $\KRh_p(T)\cong H^*(\ti\eta(\{\nai(T_i)\})).$
\end{cor}

While this may not look like an impressive simplification, it does have
a significant advantage: as a complex of modules, it is much easier to
identify trivial summands of the naive complex $\nai(T_i)$ in a way that
was not at all clear in the \mf picture.

\begin{remark}
  For instance, this allows us to replace Khovanov and Rozansky's
  exhaustive computations for invariance under Reidemeister moves with
  simple computations in Soergel bimodules (done by Rouquier
  \cite{Rou04}) for all Reidemeister moves except type I and type IIb.
  This is because we can go between any two projections by
  \begin{itemize}
  \item applying Vogel's algorithm which uses only moves of type IIb and
    passing strands through infinity on the 2-sphere (which doesn't
    change Khovanov-Rozansky, since it doesn't change the topology of
    the knot projection) to take both projections to braid-like ones
  \item applying type I moves and identities in the braid group (covered
    by Rouquier) to move from one braid projection to the other, which
    is possible by Markov's theorem.
  \end{itemize}
\end{remark}

\subsection{Braids}
\label{sec:braids}

\begin{figure}[htbp]
  \centering
    \centerline{\epsfig{figure=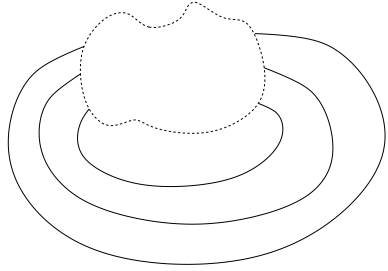, height=3.5cm}}
  \caption{The braid closure arc diagram for 3 strands, $\gamma_3$.}
  \label{fig:4}
\end{figure}

Braids have an important role to play here, especially when we wish to
consider HOMFLYPT homology, as we will in
Section~3\excise{\ref{sec:homflypt-homology}}. As we noted, one of the
best examples of a complicated acyclic tangle is a braid.  Thus,
Proposition~\ref{canop-qpm} and Theorem~\ref{acyc-qpm} show that if our
diagram $L$ is the closure of a $d$-strand braid $\sigma$ (all Seifert
circles are nested), then $\KR_p(L)$ is near-isomorphic to
$\ti\gamma_d(\nai(\sigma))$, where $\gamma_d$ is the planar arc diagram
shown in Figure~\ref{fig:4}.

In this case, the naive KR complex $\nai(\sigma)$ can be considered a
complex of bimodules over a polynomial ring $\mc S=k[x_1,\ldots,x_d]$ with
generators $x_i$ indexed by strands of our braid, and
composition of braids (which can also be written as the action of a
planar arc diagram) passes to tensor product of complexes, so
\begin{equation*}
  \nai(\sigma\sigma')\cong\nai(\sigma)\otimes_{\mc S}\nai(\sigma')
\end{equation*}
Thus, we can consider $\nai$ as a categorification of the braid group.
In fact, this is precisely the categorification of the braid group
described by Rouquier~\cite{Rou04}.  The bimodules which appear in this
complex are so-called {\bf Soergel bimodules}, which appeared in Soergel's
research on category $\mc O$.  Furthermore,
$\ti\gamma_d$ is simply the Koszul \mf $\mc{\ti Z}_{p}$ of the sequences
$\{x_i\otimes 1-1\otimes x_i\}$ and $\left\{\frac{p(x_i\otimes
    1)-p(1\otimes x_i)}{x_i\otimes 1-1\otimes x_i}\right\}$

This is an $\mathfrak{sl}_n$ version of the result of Khovanov
relating the Rouquier complex and Khovanov homology:
\begin{thm}
  The complex $H^{t}\!\left(\nai(\sigma)\otimes\mc{\ti Z}_{p}\right)$ is
  isomorphic to $\KRh_p(\bar{\sigma})$.
\end{thm}

The functor $H^t(-\otimes \ti{\mc Z}_p)$ is defined for any bimodule over
$k[x_1,\cdots x_d]$, and it would very interesting to interpret it in
terms of more familiar homological algebra. As is, there is a spectral
sequence $$H\!H^{*}(-)\Rightarrow H^t(-\otimes \ti{\mc Z}_p),$$ where
$H\!H^{*}(-)=\mathrm{Tor}^*_{\mc S\otimes \mc S^{op}}(-,\mc S)$ denotes
Hochschild homology, since $$H\!H^i(-)\cong H^i\left(-\otimes (\ti{\mc Z_p})_+\right).$$

Since there are, in all, $d!$ different indecomposable Soergel bimodules,
the complex $\nai(\sigma)$ typically has a very large number of
redundant summands.  In fact, the complex $\nai(\sigma_i)$ for a braid
generator splits after tensor product with exactly half of these
modules, which alone leads to huge number of trivial summands in
$\nai(\sigma)$ for any large braid.  This is discussed by Khovanov in the
last section of \cite{Kho05}, and will be covered in more detail in
future work by the author.

It would be even better if we could implement these cancellations
for more general acyclic tangle diagrams, since typically, a braid
representative of a given knot has many more crossings than
the smallest planar diagram of the same knot, which slows down
computation if we have to use braid diagrams.

\subsection{The IIb move.}
\label{sec:iib-move}

Another application which illustrates the power of this approach is the
IIb Reidemeister move, which creates trouble in HOMFLY homology.

\begin{figure}[htbp]
  \centering
    \centerline{\epsfig{figure=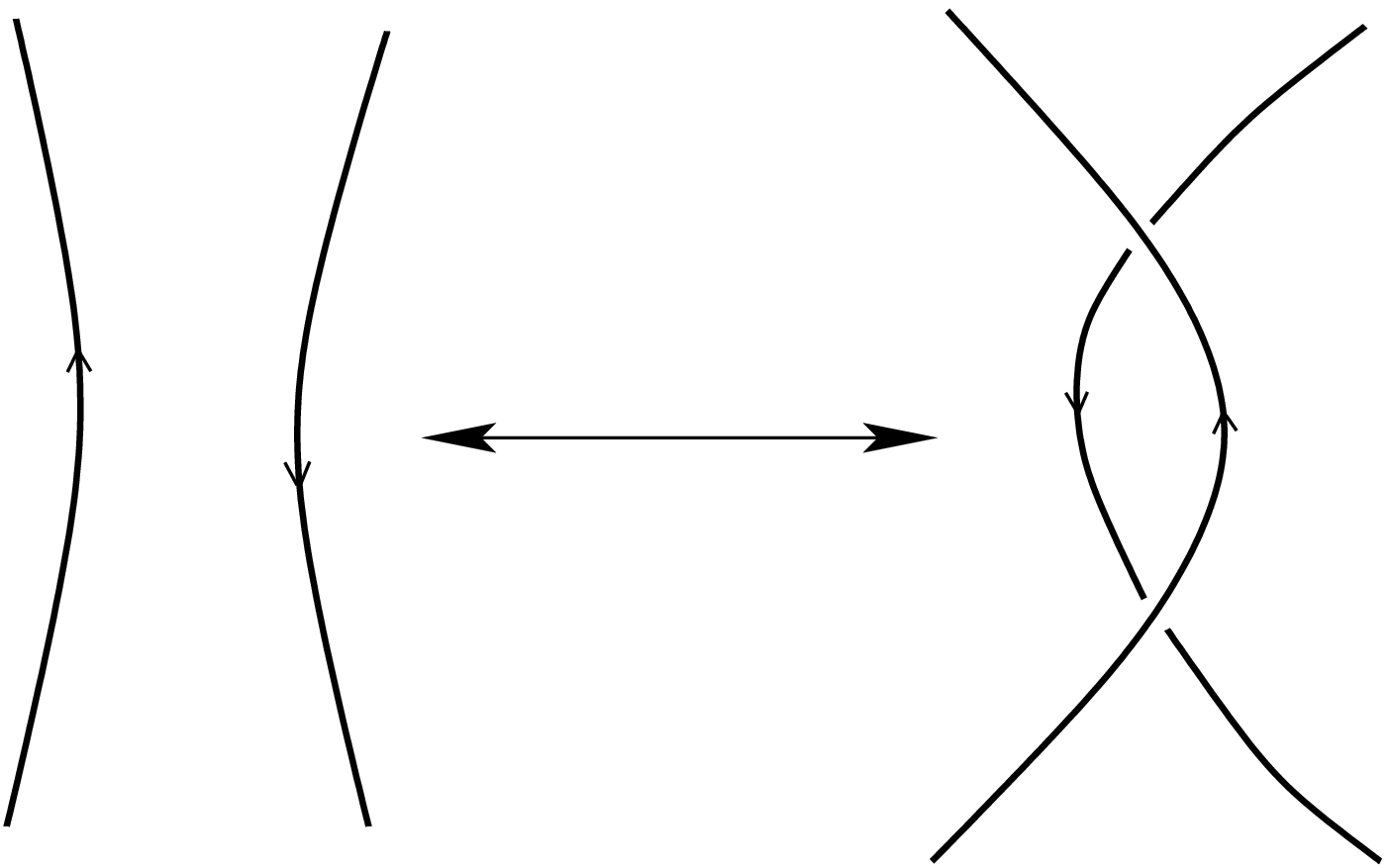, height=3cm}}
  \caption{The Reidemeister IIb move.}
  \label{fig:5}
\end{figure}

\begin{figure}[htbp]
  \centering
    \centerline{\epsfig{figure=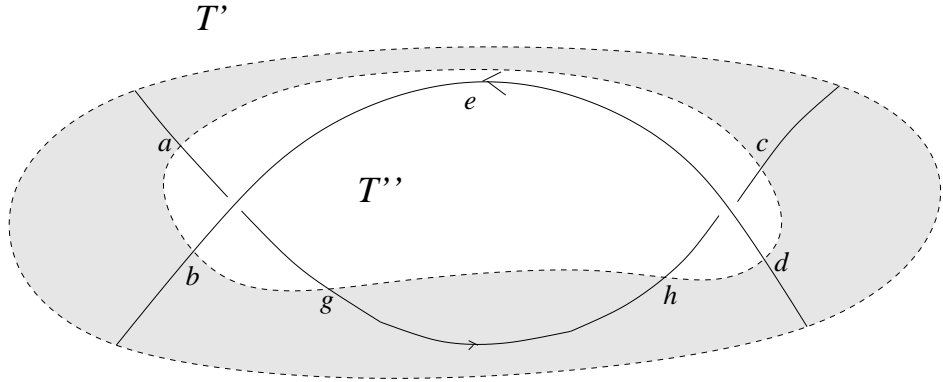, height=5cm}}
  \caption{A decomposition of the tangle $T'$.}
  \label{fig:6}
\end{figure}

We only need to understand the local picture
involving the tangles $T$ and $T'$ as shown in Figure~\ref{fig:5}.
The tangle $T$ is acyclic ,and thus quite easy to understand, but
$T'$ is not.  Thus, we will cut open the oriented cycle, consider the
naive complex of the resulting tangle $T''$, and then act with a planar
diagram $\eta$ to get $\KR_p(T)$, as shown in Figure~\ref{fig:6}.
Let $R'=k[a,b,c,d]$.  

\begin{prop}
  The complex $\KR_p(T')$ has a simplification which is a two-term
  complex of modules $\No\overset{\pi}\longrightarrow\Nt$ where
  \begin{align*}
    \No&=R'\{-2\}/(a-b+c-d,(b-d)(c-d))\\
    \Nt&=R'\{-2\}/(a-b,c-d)
  \end{align*}
  and $\pi$ is the natural projection.
\end{prop}

Note that the kernel of this surjection is isomorphic to $R'/(b-d,a-c)$,
which is, in turn near-isomorphic to $\KR_p(T)$.  Thus, if this
surjection were to split, we would have shown the invariance of KR
homology under the second Reidemeister move. However, it does not; $\No$
is indecomposable as a graded module over $R'$, since it is generated by
a single element in minimal grade.  However, as we shall see, this
calculation can still give us interesting information about this
Reidemeister move.

\begin{proof}
  Using the labels on edges above, and removing the variables $g$ and
  $h$ using the identity $a+e-g-b=c+e-h-d=0$, we get that the naive
  complex of $T''$ is the 3 term complex
\begin{equation}\label{eq:1}
\xy
  (-30,0)*{M_1}="1";
  (30,0)*{M_4}="4";
  (0,4)*{M_2}="2";
  (0,-4)*{M_3}="3";
  (0,0)*{\oplus};
  {\ar^{1} "1";"2"};
  {\ar_{e-b} "1";"3"};
  {\ar^{e-b} "2";"4"};
  {\ar_{1} "3";"4"};
\endxy
\end{equation}

where $R=k[a,b,c,d,e]$ and
\begin{align*}
  M_1&=R/(a-b,(c-d)(e-d))\\
  M_2&=R\{-2\}/(a-b,c-d)\\
  M_3&=R/((a-b)(e-b),(c-d)(e-d))\\
  M_4&=R\{-2\}/((a-b)(e-b),c-d)
\end{align*}
Since the variable $e$ corresponds to an edge that is closed, we will
only be concerned with the structure of $\nai(T'')$ as an $R'$ module,
but the purposes of calculations, it will be useful to remember the
action of $e$.

Note that $M_1$ has a decomposition as an $R'$-module into $M_1'=R'\cdot
1$ and $M_1''=R\cdot (e-d)$, and $M_4$ into $M_4'=R'\cdot 1$ and
$M_4''=R\cdot (e-b)=\mathrm{im}\, M_2$.

The module $M_3$ also has a decomposition along these lines, but a
slightly more subtle one.  We let $M_3'=R'\cdot 1+R'\cdot e$, and
$M_3''=R\cdot (e-b)(e-d)$.  Clearly, $M_3=M_3'+M_3''$, since we write
any expression with $e^n$ for $n>1$ appearing can be rewritten as the
sum of an element of $M_1''$ and a expression with a lower degree in
$e$.  On the other hand, the intersection of these submodules is
trivial, so they give a direct sum decomposition.

Thus, we can rewrite \eqref{eq:1} as
\begin{equation}
  \label{eq:2}
  \xy
  (-30,6)*{M_1'}="1";
  (-30,-6)*{M_1''}="1p";
  (0,0)*{M_2}="2";
  (0,12)*{M_3'}="3";
  (0,-12)*{M_3''}="3p";
  (30,6)*{M_4'}="4";
  (30,-6)*{M_4''}="4p";
  (-30,0)*{\oplus} ; (0,6)*{\oplus}; (0,-6)*{\oplus}; (30,0)*{\oplus};
  {\ar_{1} "1p";"3p"};
  {\ar^{\chi} "3";"4p"};
  {\ar^{e-b} "1";"3"};
  {\ar^{1} "1";"2"};
  {\ar_{e-b} "1p";"2"};
  {\ar_{1} "2";"4p"};
  {\ar_{e-b} "3p";"4p"};
  {\ar^{1} "3";"4"};
   \endxy
\end{equation}
Each module in \eqref{eq:2} is generated by a single element
over $R$ or $R'$, so most maps are induced by multiplication by a ring
element, and thus these have been denoted by the corresponding element.
The single exception is the map $\chi$, which is the natural (non-split)
projection $\chi: M_3'\to M_3'/(R'\cdot 1)$, composed with the natural
isomorphism of the target with the submodule of $M_4''$ generated over
$R'$ by $e-b$.

Now, the maps above from $M_1''$ to $M_3''$ and from $M_2$ to $M_4''$
are isomorphisms.  Let $\xi$ be the homotopy on this complex given by
the inverse of these maps, killing all other components.

One can easily calculate that $\xi\partial+\partial\xi$ is the identity
on $M_1'', M_2, M_3''$ and $M_4''$, and 0 on $M_1',M_3'$ and $M_4'$.
Thus, the complex $\nai(T'')$ is homotopic to \eqref{eq:2} with the
lower diamond removed, i.e.
\begin{equation}
  \label{eq:3}
  M_1'\overset{e-b}{\longrightarrow} M_3'\overset{1}{\longrightarrow} M_4'
\end{equation}

Let $s=a-b-c+d$.  Then acting with the planar diagram $\eta$ to connect
the ends corresponds on the matrix factorization side to tensor product
with a two-term matrix factorization, whose positive complex is just
$R'\overset{s}{\longrightarrow}R'$

Thus, we are interested in the vertical homology of
\begin{equation}\label{eq:4}
\xymatrix{M_1'\ar[r]&M_3'\ar[r]&M_4'\\
  M_1'\ar[r]\ar[u]^{s}& M_3'\ar[r]\ar[u]^{s}&M_4'\ar[u]^{s}}
\end{equation}

Since $M_1'\cong R'/(a-b)$ and $M_4'\cong R'/(c-d)$, the element $s$ is
manifestly not a zero-divisor on either of these modules.

For $M_3'$, we need only note that $M_3'$ has a basis of the form
$B_1\cup B_2$ where
\begin{equation*}
  B_1=\{a^\alpha b^\beta c^\gamma d^\delta\}\qquad B_2=\{b^\beta c^\gamma
  \cdot (e-b)\}
\end{equation*}
Note that $(c-d)(e-b)=(c-d)(b-d)$.  Thus, $s\cdot B_2$ is a linearly
independent subset of $\Span \!(B_1)$ with no $a$'s appearing, whereas
any $k$-linear combination of $s\cdot B_1$ has leading term containing
an $a$.  Thus, $s\cdot(B_1\cup B_2)$ is linearly independent, and $s$ is
not a zero-divisor, so \eqref{eq:4} has no vertical homology.

Applying Theorem~\ref{res-qpm} again, we obtain that $\KR_p(T')$ is near-isomorphic to
\begin{equation*}
  \ti M_1\overset{e-b}{\longrightarrow} \ti
  M_3\overset{1}{\longrightarrow} \ti M_4
\end{equation*}
where $\ti M_i=M_i/sM_i$.  Unlike the situation before we quotiented by
the action of $s$, the image of $\ti M_1$ in $\ti M_2$ is now
complementary to $R'\cdot 1\subset \ti M_3$.  Thus, we can do another
reducing homotopy, and see that $\KR_p(T')$ is, in fact, near-isomorphic
to the two term complex $\No=\ti M_3/\ti M_1\to \Nt=\ti M_4$ where the
map is the obvious surjection.
\end{proof}

\section{HOMFLYPT homology}
\label{sec:homflypt-homology}

Thus far, we have not discussed the triply graded theory of
Khovanov-Rozansky, defined in \cite{KR05,Kho05} which categorifies the
HOMFLYPT polynomial.  Unfortunately, this theory lacks many of the good
properties of the finite Khovanov-Rozansky theories. Most notably, it is
unclear at the moment how it changes when the IIb move is applied to the
diagram $T$, and it is not known whether or how it is functorial with
respect to embedded cobordisms.

Typically, because of the issues surrounding the IIb move, HOMFLYPT
homology is defined only using a braid representation of the knot.  We
would rather take the perspective that it is an invariant of tangle
diagrams with good properties under Reidemeister moves.  We will briefly
discuss how our computational schema extends to HOMFLYPT homology, and
use this to obtain some information about how this homology theory
reacts to the IIb move.

Having already built up our machinery around $\KR_p$, defining HOMFLYPT
homology is simple: we consider $\KR_0$, that is, the bicomplex given by
considering the positive complex of each matrix factorization in
$\KR_p$.  As before, we take homology first in the ``matrix
factorization'' direction, and then take homology of the resulting chain
complexes, and apply a grading shift of $-w+b$ in the polynomial
grading, $w+b-1$ in the ``\mf'' grading and $w-b+1$ in the
``cohomological'' grading.

We will draw our complexes with the ``matrix factorization'' direction
being horizontal and the ``cohomological'' direction being vertical, and
henceforth use these terms to describe them.  We denote the $i$th
``horizontal'' homology (remember this is itself a single chain complex)
$\Hv iC$ and denote ``horizontal then vertical'' homology by
$\Hhv{i,j}{C}=\Hh j{\Hv i{C}}$ for any bicomplex $C$.

In line with our previous notation we denote $\Hv{i}{\KR_0(T)}$ by
$\KRh_0^i(T)$.  The direct sum of these over $i$ is the (now bi-graded)
complex $\KRh_0(T)$

Since each group $\Hhv{i,j}{\KR_0(T)}$ is still a graded module over
$\mc S_T$, we have three gradings, one inherited from the polynomial ring
$\mc S_T$, and two cohomological gradings.

Most of our previous theorems remain true, and in fact, are much easier
to prove, since we no longer need to consider matrix factorizations.  We
have lost \excise{both functoriality and} invariance under Reidemeister moves,
but if we consider $\KR_0$ only as an invariant of tangle diagrams,
essentially everything works as before.
\begin{thm}\label{homf-can}
  For a sequence of tangles $\{T_i\}$ and compatible oriented
  planar arc diagram $\eta$, we have 
  \begin{equation*}
    \KR_0(\ti\eta(\{T_i\}))\cong \ti\eta(\{\KR_0(T_i)\})
  \end{equation*}

  If $\KR_0(T_i)\to\nai(T_i)$ is a quasi-isomorphism, term-wise, then
  $\KR_0(\ti\eta(\{T_i\}))\to\ti\eta(\nai(T_i))$ is also a
  quasi-isomorphism. 
\end{thm}

In fact, since complexes are easier to deal with than matrix
factorizations, in this case, our invariants can be understood in terms
of standard homological algebra.

For instance, if we factor a closed diagram $T$ into a planar arc
diagram $\eta$ acting on a set $\{T_i\}$ of acyclic tangles (as shown
in Figure~\ref{fig:2}), then we can partition the set $X$ of endpoints
of the tangles $T_i$ into pairs $\al_+,\al_-$ where $\al$ ranges over
$\mc A(\eta)$, the set of arcs of $\eta$.  Let $\mc
S=\bigotimes_{i}\mc S_{T_i}=k[X]$ be the ring over which the bicomplex
$\bigotimes_i\nai(T_i)$ is defined, and $\mc U=k[\mc A(\eta)]$. In this case, $\mc S$ can be written as a tensor product
$\mc U\otimes \mc U^{op}$ with the left
action of $t_{e}$ being $t_{e_+}$ and its right action being
$t_{e_-}$.  Thus, we can also consider $\bigotimes_i\nai(T_i)$ as a
complex of bimodules over $\mc U$.

\begin{prop} The horizontal homology complex 
  $\KRh_0^i(T)$ is naturally isomorphic to
  the Hoch\-schild homology
  $H\!H_i^{\mc U}\left(\bigotimes_i\nai(T_i)\right)$ where $H\!H_i^{\mc U}(-)$
  is applied term-wise.
\end{prop}
\begin{proof}
  By Theorem~\ref{homf-can}, $\KRh_0^i(T)\cong
  \Hv{i}{\ti\eta\left(\{\nai(T_i)\}\right)}$. The bicomplex
  $\ti\eta(\{\nai(T_i)\})$ is simply the tensor product of
  $\otimes_i\nai(T_i)$ (thought of as a vertical bicomplex) with the
  horizontal bicomplex which is the Koszul complex of $\left(t_e\otimes
    1-1\otimes t_e\right)$ over $\mc S$.  This complex is a free
  resolution of $\mc U$ as a module over $\mc S$ (i.e. as a bimodules
  over itself), so the homology of the tensor product of this complex
  with a bimodule over $\mc U$ is precisely the Hochschild homology of
  that bimodule.
\end{proof}
Note that if $T_1$ is a braid, and $T$ its closure, this theorem
reduces precisely to \cite[Theorem 1]{Kho05}.
  
As this result suggests, we can use a stronger notion of equivalence in
this HOMFLY case.  Recall that $D^b(\mc S_T)$, the derived category of
$\mc S_T$-modules, is the category of complexes of $\mc S_T$-modules, with
a formal inverse to each quasi-isomorphism added.  Since this category
is additive, the homotopy category $\hd{\mc S_T}$ of complexes in
$D^b(\mc S_T)$
is well defined.

Of course, the operation of $\Hv{*}{-}$ is still well defined over
$\hd{\mc S_T}$, and results in a series of chain complexes of $\mc S_T$-modules.
Furthermore, any nullhomotopic map in $\hd{\mc S_T}$ induces nullhomotopic
maps on these complexes.  Thus $\Hhv{*,*}{-}:\hd{\mc S_T}\to \mc S_T-\mathsf{mod}$
is a well-defined functor.  Furthermore, using
Proposition~\ref{T:tens-qi}, we see that tensor product with a bicomplex
of $\mc S_T\otimes \mc S_{T'}$-modules which is projective as a $\mc S_T$-module $D$
defines a functor $D\otimes -:\hd{\mc S_T}\to\hd{\mc S_T'}$.

Thus, we have the proposition
\begin{prop}
  The bigraded complex $\KRh_0(T)$ only depends (up to homotopy) on the
  class of $\KR_0(T)$ in $\hd{\mc S_T}$.  Furthermore, the canopolis
  $\EuScript{M}_0$ defined in Section~\ref{sec:canopo-whats} descends to
  a canopolis structure $\EuScript{M}_0'$ on the categories $\hd{\mc S_T}$.
\end{prop}

To show the power of this approach, let us give a proof of the
invariance of HOMFLYPT using it.  Recall the famous theorem of Markov:
\begin{thm}
  \emph{(Markov, \cite{Bir74})} Two closed braid projections represent
  the same knot if and only if they are related by isotopy, identities
  in the braid group, and type I Reidemeister moves.
\end{thm}

\begin{figure}[htbp]
  \centering
    \centerline{\epsfig{figure=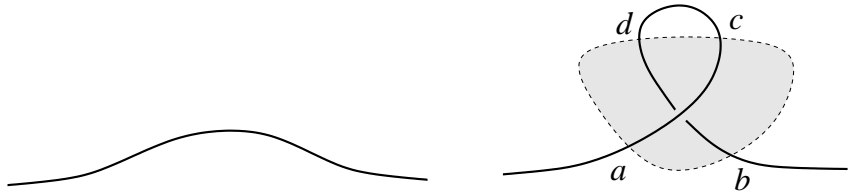, height=2.5cm}}
  \caption{The Reidemeister I move.}
  \label{fig:7}
\end{figure}

\begin{proof}[Proof of invariance of HOMFLYPT homology]
  Isotopy does not affect the structure of the KR complex, and relations
  in the braid group have been dealt with by Rouquier \cite{Rou04}.
  Thus, we need only consider type I moves.  In both cases, we will only
  obtain invariance with respect to a grading shift.  This is accounted
  for in the global grading shift, which also depends on the diagram.

  The left side of a type I move is simply the bicomplex which is
  $k[a,b]/(a-b)$ in cohomological grading $(0,0)$.

  The right side is the tensor product of two 2-term bicomplexes, one
  being the naive complex of the crossing, and the other corresponding
  to the diagram closing one end of the crossing.  Let
  \begin{align*}
    \mc S&=k[a,b,c,d]\\
    M_0&=\mc S/(a-b,c-d)\\
    M_1&=\mc S/\left(a+c-b-d,(a-b)(a-c)\right)
  \end{align*}
  Note that multiplication by $a-c$ defines a map $\rho:M_0\{2\}\to M_1$,
  which is injective.  Let $M_1''$ denote the image of $\rho$.  As a
  $k[a,b]$-submodule, we have a decomposition of the form
  $M_1''=\oplus_{i=0}^\infty k[a,b]\cdot c^i(a-c)$, and this image has a
  complement $M_1'=k[a,b]\cdot 1$ i.e.  $M_1\cong M_1'\oplus M_1''$ (of
  course, $M_1'$ is not a $k[a,b,c]$-submodule).

  For the negative move, the bicomplex $\KR_0(T_+)$ is
  \begin{equation*}
     \xy
     (-10,7)*{M_0}="1";
     (-10,-5)*{M_1''\{-2\}}="2"; 
     (-7,-8.5)*{\oplus}; 
     (-4,-12)*{M_1'\{-2\}}="2b";
     (15,7)*{M_0}="3";
       (15,-5)*{M_1''\{-2\}}="4"; 
       (18,-8.5)*{\oplus}; 
       (21,-12)*{M_1'\{-2\}}="4b";
     {\ar^{\rho} "1";"2"};
      {\ar_{\rho} "3";"4"};
      {\ar_{a-b} "2b";"4b"};
       \endxy
  \end{equation*}
  Since the vertical map induces an isomorphism from $M_0$ to $M_1'$,
  after removing the null-homotopic summand, $\KR_0(T_+)$ is just
  the horizontal complex $M_1'\overset{a-b}{\longrightarrow}M_1'$, which
  is, in turn, quasi-isomorphic to $k[a,b]/(a-b)\{-2\}$ with a vertical
  shift of $1$.
  
  For the positive move, we must be a bit more subtle. We start with the
  complex
   \begin{equation*}
     \xy
     (-15,12)*{M_1'}="1a"; (-20,8.5)*{\oplus}; (-25,5)*{M_1''}="1b";
      (-20,-5)*{M_0}="2"; {\ar_{a-c} "1b";"2"}; {\ar^{1} "1a";"2"}; 
      (25,12)*{M_1'}="3a"; (20,8.5)*{\oplus}; (15,5)*{M_1''}="3b";
      (20,-5)*{M_0}="4"; {\ar_{a-c} "3b";"4"}; {\ar^{1} "3a";"4"};
      {\ar^{a-b} "1a";"3a"};
     \endxy
  \end{equation*}
  Note that by the decomposition mentioned earlier
  $M_0\cong M_0\{2\}\oplus k[a,b]/(a-b)$.  In the derived category, we
  can replace $k[a,b]/(a-b)$ by the complex
  $M_1'\overset{a-b}{\longrightarrow}M_1'$.  Thus, we can write another
  representative of this complex in the derived category which is
  \begin{equation*}
      \xy
      (-15,12)*{M_1'}="1a"; (-20,8.5)*{\oplus}; (-25,5)*{M_1''}="1b";
      (-15,-5)*{M_1'}="2a"; (-20,-8)*{\oplus};
      (-25,-12)*{M_0}="2b";
      (15,12)*{M_1'}="3a";(20,8.5)*{\oplus}; (25,5)*{M_1''}="3b";
      (15,-5)*{M_1'}="4a";(20,-8)*{\oplus}; (25,-12)*{M_0\{2\}}="4b";
      {\ar^1 "1a";"2a"}; {\ar_{a-c} "1b";"2b"}; {\ar_{1} "1a";"2b"};
      {\ar^1 "3a";"4a"}; {\ar_{1} "3b";"4b"};
      {\ar_{a-b} "2a";"4a"};  {\ar^{a-b} "1a";"3a"}
    \endxy
  \end{equation*}
  The the top row of $M_1'$'s and its image in the bottom row form a
  null-homotopic subcomplex, as does the far right column.  Canceling
  these off, we obtain a split injection (by the same decomposition we
  used before), with cokernel $k[a,b]/(a-b)$.  Thus, $\KR_0(T_-)$ is
  equivalent to a single copy of $k[a,b]/(a-b)$ but now with a
  horizontal shift of $-1$.  
\end{proof}

\subsection{IIb Again}
\label{sec:iib-again}

Let us return to the IIb move.  As we mentioned earlier, the behavior of
HOMFLYPT homology under this Reidemeister move is not well understood.
Using the results of the previous sections, we will make some headway
toward understanding HOMFLYPT homology for general diagrams.

As we showed in Section~\ref{sec:iib-move}, the homology $\KR_0(T')$ is
near-isomorphic (or more precisely, derived homotopic) to a two term
complex $\No\overset{\pi}{\longrightarrow}\Nt$, which is
quasi-isomorphic but not homotopic to $\KR_0(T)$.

Thus, if $K$ is any link diagram, with $T$ a subdiagram which is isotopic
to that on the left side of the IIb move, and $K'$ the diagram which
results after a IIb move, we can construct representatives in $\hd{\mc S_K}$
of $\KR_0(K)$ and $\KR_0(K')$ such that there is an injective map
$\iota:\KR_0(K)\to\KR_0(K')$, with the cokernel of $\iota$ given by the
tensor product
$\KR_0(K\backslash T)\otimes
\Big(\Nt\overset{\mathrm{id}}{\longrightarrow}\Nt\Big)$.
\nc{\coker}{\mathrm{coker}\,}

Of course, we are interested in understanding $\ker\left(\Hhv
  {i,j}{\iota}\right)$ and $\coker\left(\Hhv {i,j}{\iota}\right)$. 

\begin{prop}
  There is a spectral sequence $E_n^{i,j}$ such that
  \begin{equation*}
    \label{eq:5}
    E^{i,j}_{2}=
    \begin{cases}
      \ker\left(\Hhv
  {i,j}{\iota}\right) & i=3k\\
    \coker\left(\Hhv {i,j}{\iota}\right) & i=3k+1\\
    0 & i=3k+2
    \end{cases}
    \Rightarrow E_{\infty}^{i,j}=0
  \end{equation*}

  Furthermore, we have a complex 
  \begin{equation*}
    \cdots\longrightarrow
    \Hhv{i,j}{K}\overset{\Hhv{i,j}\iota}{\longrightarrow}
    \Hhv{i,j}{K'}\overset{\al_i}{\longrightarrow}
    \Hhv{i-1,j-1}{K}\overset{\Hhv{i-1,j-1}\iota}
    \longrightarrow\Hhv{i-1,j-1}{K'} \longrightarrow\cdots
  \end{equation*}
  which exact if and only if the spectral sequence above collapses at
  the $E_3$-term (i.e. $d^i=0$ for $i\geq 3$).
\end{prop}

Thus, if this spectral sequence collapses early on, the IIb move does
relatively little damage; part of the HOMFLY homology shifts by one in
horizontal and vertical grading, and part of it does not.
Unfortunately, as of the moment, we have not able to obtain any real control
over the higher differentials.

In the author's view, the most optimistic hope is that the spectral
sequence does collapse at $E_3$, and that $\Hhv{i,j}\iota$ is a
isomorphism if the crossing strands of the IIb move lie on the different
Seifert circles and $0$ (and thus $\al_i$ is an isomorphism) if they lie
on the same Seifert circle.  This would imply that with the grading
shifts described earlier, HOMFLY homology was independent of the diagram
chosen.

Another weaker possibility is that the spectral sequence collapses, and
there is some good description of the differentials less clean than the
hope above. 

Weaker still is the hope that $d^{3n}=0$ for all $n\geq 1$.  This would
at least imply that the total rank does not change under the IIb move,
and in fact that each piece of the homology shifts by a vector lying in
one of two affine rays in $\Z^3$. 
\begin{proof}
  By the usual yoga, since we have a short exact sequence of
  bi-complexes (though of a vertical complexes of horizontal chain
  complexes) we obtain a long exact sequence of vertical complexes
\begin{equation*}
  \cdots\to\Hv i{\KR_0(K)}\to\Hv i{\KR_0(K')}\to\Hv i{\coker\iota}\to\Hv {i-1}{\KR_0(K)}\to\cdots
\end{equation*}

Now, think of this long exact sequence itself as a bicomplex.  Then we
have a pair of spectral sequences converging from ``vertical then
horizontal'' homology to total homology and from ``horizontal then
vertical'' homology to total homology.  Since this is an exact sequence
of complexes, the horizontal homology, and thus total homology is
trivial.  On the other hand, taking vertical homology first, we obtain
a spectral sequence converging to 0 with $E^1$ page
\begin{equation*}
  \xymatrix@R=.4cm{ 
&\vdots&\vdots&\vdots&\vdots&\\
\cdots\ar[r] & \Hhv {i,j}{\KR_0(K)}\ar[r]^{\Hhv {i,j}{\iota}}&\Hhv {i,j}{\KR_0(K')}\ar[r]&0\ar[r]&\Hhv
  {i-1,j}{\KR_0(K)} \ar[r]& \cdots \\
 \cdots\ar[r] &\ar[r]^{\Hhv {i,j-1}{\iota}} \Hhv {i,j-1}{\KR_0(K)}&\ar[r]\Hhv {i,j-1}{\KR_0(K')}&\ar[r]0&\ar[r]\Hhv
  {i-1,j-1}{\KR_0(K)} & \cdots \\
&\vdots&\vdots&\vdots&\vdots&}
\end{equation*}

Taking homology, we see that the $E^2$ page of the same spectral
sequence is
\begin{equation*}
   \xymatrix@R=.4cm{ 
&\vdots\ar[drr]&\vdots\ar[drr]&\vdots&\vdots&\\
\cdots \ar[drr]& \ker \Hhv {i,j}{\iota} \ar[drr]& \coker \Hhv {i,j}{\iota} \ar[drr] & 0\ar[drr] & \ker \Hhv {i-1,j}{\iota}& \cdots\\
\cdots & \ker \Hhv {i,j-1}{\iota}\ar[drr] & \coker \Hhv {i,j-1}{\iota} \ar[drr]& 0 & \ker \Hhv {i-1,j-1}{\iota}& \cdots\\
&\vdots&\vdots&\vdots&\vdots&}
\end{equation*}
Thus, this the desired spectral sequence.  Furthermore, $d^2$ defines a
map $\al_i:\Hhv {i,j}{\KR_0(K')}\to \Hhv {i-1,j-1}{\KR_0(K)}$ by
composition with the projection to $\coker \Hhv {i,j}{\iota}$ and the
inclusion of $\ker \Hhv {i-1,j-1}{\iota}$.  These maps define a complex
by definition, and this complex is exact if and only if $d^2$ is an
isomorphism.  Since this spectral sequence converges to 0, it collapses at
$E^n$ if and only if $d^{n-1}$ is an isomorphism.
\end{proof}
\bibliographystyle{../halpha} 
\bibliography{../gen}

\def\cprime{$'$}
\begin{thebibliography}{Str06b}

\bibitem[Bir74]{Bir74}
Joan~S. Birman.
\newblock {\em Braids, links, and mapping class groups}.
\newblock Princeton University Press, Princeton, N.J., 1974.
\newblock Annals of Mathematics Studies, No. 82.

\bibitem[BN02]{BarN02}
Dror Bar-Natan.
\newblock On {K}hovanov's categorification of the {J}ones polynomial.
\newblock {\em Algebr. Geom. Topol.}, 2:337--370 (electronic), 2002.

\bibitem[BN05]{BarN05}
Dror Bar-Natan.
\newblock Khovanov's homology for tangles and cobordisms.
\newblock {\em Geom. Topol.}, 9:1443--1499 (electronic), 2005.

\bibitem[Eis80]{DE80}
David Eisenbud.
\newblock Homological algebra on a complete intersection, with an application
  to group representations.
\newblock {\em Trans. Amer. Math. Soc.}, 260(1):35--64, 1980.

\bibitem[Eis95]{Ei97}
David Eisenbud.
\newblock {\em Commutative algebra}.
\newblock Springer-Verlag, New York, 1995.
\newblock With a view toward algebraic geometry.

\bibitem[FKS]{FKS05}
Igor Frenkel, Mikhail Khovanov, and Catharina Stroppel.
\newblock {A categorification of finite-dimensional irreducible representations
  of quantum sl(2) and their tensor products}, arXiv:math.QA/0511467.

\bibitem[Jon99]{Jon99}
Vaughan F.~R. Jones.
\newblock {Planar algebras, I}, 1999, arXiv:math.QA/9909027.

\bibitem[Kho05]{Kho05}
Mikhail Khovanov.
\newblock {Triply-graded link homology and Hochschild homology of Soergel
  bimodules}, 2005, arXiv:math.GT/0510265.

\bibitem[KR04]{KR04}
Mikhail Khovanov and Lev Rozansky.
\newblock {Matrix factorizations and link homology}, 2004,
  arXiv:math.QA/0401268.

\bibitem[KR05]{KR05}
Mikhail Khovanov and Lev Rozansky.
\newblock {Matrix factorizations and link homology II}, 2005,
  arXiv:math.QA/0505056.

\bibitem[Ras]{Ras06}
Jacob Rasmussen.
\newblock {Some differentials on Khovanov-Rozansky homology},
  arXiv:math.GT/0607544.

\bibitem[Rou04]{Rou04}
Raphael Rouquier.
\newblock {Categorification of the braid groups}, 2004, arXiv:math.RT/0409593.

\bibitem[Soe92]{Soe92}
Wolfgang Soergel.
\newblock The combinatorics of {H}arish-{C}handra bimodules.
\newblock {\em J. Reine Angew. Math.}, 429:49--74, 1992.

\bibitem[Str06a]{Str06b}
Catharina Stroppel.
\newblock {Perverse sheaves on Grassmannians, Springer fibres and Khovanov
  homology}, 2006, arXiv:math.RT/0608234.

\bibitem[Str06b]{Str06a}
Catharina Stroppel.
\newblock {TQFT with corners and tilting functors in the Kac-Moody case}, 2006,
  arXiv:math.RT/0605103.

\end{thebibliography}
\end{document}